\documentclass[12pt,a4paper]{article} 

\usepackage[utf8]{inputenc}
\usepackage[T1]{fontenc}
\usepackage{amsmath,amssymb,amsthm} 
\usepackage{graphicx}
\usepackage{geometry}
\usepackage{hyperref}
\usepackage{color}

\newtheorem{theorem}{Theorem}[section]

\newtheorem{definition}[theorem]{Definition}

\theoremstyle{definition}

\title{Graph Neural Networks for Community Detection in Graph Signal Analysis}
\author{Roberto Cavoretto, Alessandra De Rossi, Enrico Montini\thanks{Department of Mathematics \lq\lq Giuseppe Peano\rq\rq, University of Torino, via Carlo Alberto 10, 10123 Torino, Italy}}
\date{\texttt{roberto.cavoretto@unito.it, alessandra.derossi@unito.it, enrico.montini2000@gmail.com}}

\begin{document}

\maketitle

\begin{abstract}
Community detection is a central problem in graph analysis, with applications ranging from network science to graph signal processing. In recent years, Graph Neural Networks (GNNs) have emerged as effective tools for learning low-dimensional representations of graph-structured data and have shown strong performance in clustering tasks, particularly on large and high-dimensional graphs. This paper investigates the use of GNN-based community detection within a graph signal interpolation framework. After reviewing the main classes of GNN architectures for community detection according to a standard taxonomy, we integrate the resulting graph communities into a Partition of Unity Method (PUM) for interpolation with Graph Basis Functions (GBFs). In this approach, GNN-derived communities are used to construct local subdomains on which GBF interpolants are computed and subsequently combined into a global approximation. Numerical experiments on benchmark 
urban network examples demonstrate that the proposed combination of GNN-based clustering and GBF-PUM interpolation yields accurate signal reconstructions. The results indicate that deep learning-based community detection can provide effective graph partitions for localized interpolation schemes, supporting its use in scalable graph signal analysis.
\end{abstract}

\vspace{0.5cm} 

\section{Graphs and Graph Neural Networks}

\subsection{Basics of Graph Theory}

\begin{definition} A graph $G$ is defined by an ordered couple $(V(G), E(G))$, where $V(G)$ and $E(G)$ are, respectively, the set of $n$ nodes and $m$ edges of the graph, such that every edge is connected to a couple of nodes \cite{bondy}.
\end{definition}

In an indirect graph, the couple of nodes that identify an edge is not ordered; an indirect graph without loops (edges where the extremes are coincident) or parallel edges (edges with the same pair of nodes as extremes) is called simple. For the rest of this paper, we will only consider simple graphs. 

\begin{definition} Let $G = (V(G), E(G))$ be a graph. The $n \times n$ adjacency matrix $\textbf{A}$ and degree matrix $\textbf{D}$ associated with the graph are respectively defined as
\begin{equation*} \textbf{A}_{ij} = \begin{cases} 0, & \mbox{if } \{v_i, v_j\} \notin E(G), \\ 1, & \mbox{if } \{v_i, v_j\} \in E(G), \end{cases}
\hspace{1cm}
\textbf{D}_{ij} = \begin{cases} 0, & \mbox{if } i \neq j, \\ d_j = \sum_{k = 1}^n \textbf{A}_{jk}, & \mbox{if } i = j. \end{cases}
\end{equation*}
From these matrices we can obtain the Laplacian matrix $\textbf{L}$ and the normalized Laplacian matrix $\textbf{L}_n$ associated with graph $G$:
\begin{equation*}\textbf{L} = \textbf{D} - \textbf{W}, \hspace{1cm}\textbf{L}_n = \textbf{I}_n - \textbf{D}^{-\frac{1}{2}} \textbf{A} \textbf{D}^{-\frac{1}{2}}.
\end{equation*}
Here, $\textbf{I}_n$ is the $n \times n$ identity matrix.
\end{definition}

\subsection{Community Detection in Graphs}

The distribution of edges in graphs can be extremely inhomogeneous \cite{fortunato}: some regions may exhibit dense connectivity, while others may be relatively sparse. The task of community detection, or clustering, requires one to divide the graph into a set of communities $\overline{C} = \{C_1, \dots, C_J\}$ that completely cover the graph, so that the density of edges in a community is higher than in the entire graph. If the communities do not overlap, we call $\overline{C}$ a partition, while if they do we refer to it as a cover. 

A quality function assigns a value to every partition of the graph based on its properties. The most utilized quality function is the modularity (MOD) introduced in \cite{newman} and defined as
\begin{equation*} \label{modularity} Q(\overline{C}) = \frac{1}{2m} \sum_{v_i, v_j \in V(G)} (\textbf{A}_{ij} - \textbf{P}_{ij}) \delta(C_i, C_j),
\hspace{1cm} \textbf{P}_{ij} = \frac{d_i d_j}{2m}.
\end{equation*}


Community detection has attracted considerable attention in recent years, leading to the development of a wide range of algorithmic approaches \cite{fortunato}. In this paper, we focus on methods based on Graph Neural Networks (GNNs).

\section{GNNs for Community Detection}

\subsection{Introduction to GNNs}
The irregular structure and often large scale of graph data have posed significant challenges for conventional deep learning methods \cite{wu}. Therefore, new classes of neural networks customized to work on graphs have been developed, grouped under the common name of GNNs.

Given a graph $G$ and an initial node representation matrix $\textbf{X} \in \mathbb{R}^{n \times C}$, the $k$-th layer of the GNN generates a new representation $\textbf{H}^k_v$ of the node $v$ by aggregating the previous representations of all $v$-adjacent nodes from the $(k-1)$-th layer and combining it with the representation of $v$.
By optimizing a suitable cost function, the model learns its parameters and improves the representation.

GNNs are able to learn an effective low-dimensional representation of graph data, which can be useful in numerous fields \cite{wu}. For community detection in particular, the ability of GNNs to learn lower-dimensional embeddings from high-dimensional graph data led deep learning methods to obtain much more accurate clustering results when compared with more classic methods \cite{jin,su}.

\subsection{GNN Methods for Community Detection on Graphs}
To present GNN-based clustering methods, we organize them according to a common taxonomy drawn from existing surveys on the topic \cite{jin}, \cite{su}.

\subsubsection{Convolutional GNNs}


In convolutional GNNs, node representations are aggregated through a convolution operation defined by an appropriate graph filter.

Convolutional GNNs are mostly used in a semi-supervised way, but by combining them with other clustering techniques, it is possible to operate them in an unsupervised manner; since community detection is an unsupervised task, especially when real communities are not known, we focus primarily on unsupervised methods.

A way to generate a good partition is by using a low-pass filter, which gives a smoother signal after the convolution. As an example, the AGC method \cite{zhang} finds iteratively 
an appropriate function to construct an effective graph filter.

The most common convolutional structure is GCN \cite{kipf}. This method overlaps many convolutional layers of the form
\begin{equation*}\textbf{H}^{(k + 1)} = \sigma(\hat{\textbf{A}} \textbf{H}^{(k)} \boldsymbol{\Theta}), \hspace{1cm}
\hat{\textbf{A}} = (\textbf{D} + \textbf{I}_n)^{-\frac{1}{2}}(\textbf{A} + \textbf{I}_n)(\textbf{D} + \textbf{I}_n)^{-\frac{1}{2}}.
\end{equation*}
Several clustering methods are based on the GCN structure. 
For example, both DGCluster \cite{bhowmick} and NOCD \cite{shchur} methods generate communities from the GCN-generated representations by approximating modularity or creating a similarity matrix between nodes. The NOCD method is also capable of identifying overlapping community structures.

\subsubsection{Graph attention networks}

In a graph attention network (GAT) \cite{velickovic}, node features are aggregated via trainable attention weights: the $k$-th GAT layer takes the form 
\begin{equation*}\textbf{H}^{(k)} = \sigma\left(\mathbf{\Lambda}^{(k)} \mathbf{\Theta}^{(k)} \textbf{H}^{{(k-1)}^T}\right),
\end{equation*}
where element $\lambda^{(k)}_{ij}$ of the matrix $\mathbf{\Lambda}^{(k)}$ is the attention coefficient between nodes $v_i$ and $v_j$. The attention mechanism proves to be very useful in community detection, as it can highlight connections between nodes in the same community while punishing those between separate communities \cite{li}. For example, the AOCD method \cite{sismanis}, which builds on the previously discussed NOCD framework, replaces the GCN component with a $q$-headed GAT, consisting of $q$ independent attention mechanisms whose outputs are subsequently averaged.

GAT is mainly used in community detection on multiplex networks, where edges belong in different classes. Here, the attention mechanism can be used to distinguish between different types of connections \cite{su}.

\subsubsection{Auto-Encoders}

Auto-Encoders (AE) are unsupervised networks formed by an encoder, that generates a representation matrix $\textbf{Z} \in \mathbb{R}^{n \times F}$ of the graph, and a decoder, that recreates a new adjacency matrix $\overline{\textbf{A}}$ or a characteristic matrix $\overline{\textbf{X}}$. The model learns its parameters via a comparison between the real and the reconstructed matrices.

The AEs for community detection can be divided in four main classes:
\begin{itemize}
    \item Stacked AEs, as the Semi-DNR method \cite{yang}, stack together multiple AE layers, gradually increasing the quality of the representation;
    
    \item Sparse AEs, as the GraphEncoder method \cite{tian}, add a penalty cost $\Omega(\textbf{Z})$ to the cost function to maintain a sparse output, easier to store and work with;
    
    \item Denoising AEs, as the MGAE method \cite{wang}, train the encoder on a noisy input, so to get a representation which is more robust to noise;
    
    \item Variational AEs, as the MAVGAE method \cite{salha}, use the popular VGAE method \cite{kipf_2}, an AE with probabilistic encoder and decoder, adapting it to clustering.
\end{itemize}

AEs can also be used to improve representation generated by other methods, as an example by turning a semi-supervised method into a non supervised one: this is the case with the CDBNE method \cite{zhou}, which uses a GAT as the encoder. 

\subsubsection{Graph Adversarial Networks}

A Graph Adversarial Network (GAN) \cite{goodfellow} works on the graph in a usually unsupervised way, learning a distribution $p_{data}$ of the information with the combination of a generator $G(z, \theta_G)$ that generates fake data from a seed $z$, usually sampled from a normal distribution, and a discriminator $D(x, \theta_D)$ which is trained to distinguish between the fake data generated via $G$ and the real graph data. Parameters $\theta_G$, $\theta_D$ are learned via a minimax cost function, where $G$ is trained to deceive the discriminator while $D$ tries to identify real data. 

In community detection, GANs can be used in collaboration with other methods to refine their representation, such as the AMIL method \cite{he_2}, which uses a GAN to compete with a Variational Auto-Encoder, improving the representation generated from the first one.  

\section{Signal Interpolation on Graphs}

\begin{definition}
Given a graph $G$, a signal on $G$ is a function $x: V(G) \longrightarrow \mathbb{R}$. This can also be interpreted as a vector $\textbf{x} \in \mathbb{R}^n$. The space of all signals on $G$ is defined as $L(G)$.
\end{definition}

On $L(G)$ we can define the Fourier transform and the convolution operation between signals respectively as 
\begin{equation*}\hat{\textbf{x}} = \textbf{U}^T \textbf{x},
\hspace{1cm} \textbf{x} * \textbf{y} = \textbf{C}_{\textbf{x}} \textbf{y} = \textbf{U} \textbf{M}_{\hat{\textbf{x}}} \textbf{U}^T \hat{\textbf{y}},
\end{equation*}
where $\textbf{U}$ is the eigenvector matrix of the normalized Laplacian $\mathbf{L}_n$. 

\begin{definition} Given a signal $x \in L(G)$ and a set of nodes $W = \{w_1, \dots, w_N\} \subset V(G)$, an interpolation of $x$ on $W$ is a signal $I_Wx \in L(G)$ which satisfies the interpolation conditions
\begin{equation*}I_Wx(w_i) = x(w_i), \hspace{0.5cm} i = 1, \dots, N.
\end{equation*}
\end{definition}

A way to solve the interpolation problem is to use a positive definite Graph Basis Function (GBF), namely a signal $f \in L(G)$ such that $\hat{f}_k > 0$, $\forall k \in \{1, \dots, n\}$ \cite{erb}. A GBF gives a unique interpolation $I_Wx$ as
\begin{equation*}I_Wx(v) = \sum_{i = 1}^N c_i \textbf{C}_{e_{i}}f(v), 
\hspace{1cm} 
e_{i}(v) = \delta_{w_{i}}(v) = \begin{cases} 0, & \mbox{if } v \neq w_i \\ 1, & \mbox{if } v = w_i\end{cases} \hspace{0.2cm} \in L(G).
\end{equation*}
The interpolation conditions can be formulated as an $N \times N$ linear system, determining coefficients $c_i$.
Since solving such a system can be computationally costly for large $N$, a possible approach is to use a Partition of Unity Method (PUM) on the graph \cite{cavoretto,cavoretto_2}. Starting from a partition $\overline{C}$, the communities are enlarged via a factor $R \geq 0$, obtaining a cover $\overline{V} = \{V_1, \dots, V_J\}$; on this cover, we construct a partition of unity $\{\phi^{(j)}\}_{j = 1}^J$, i.e. a family of functions $\phi^{(j)} \in L(G)$ such that
\begin{equation*}\text{supp}(\phi^{(j)}) \subset V_j, \hspace{1cm}\phi^{(j)} \geq 0, \hspace{1cm} \sum_{j = 1}^J \phi^{(j)}(v) = 1, \hspace{0.3 cm} \forall v \in V(G).
\end{equation*}
The global interpolation is then determined as
\begin{equation*} x_*(v) = \sum_{j = 1}^J \phi^{(j)}(v) x_*^{(j)}(v), \hspace{0.5cm} \forall v \in V(G),
\end{equation*}
where $x_*^{(j)}$ is the local interpolation on $V_j$. As the GBF-PUM algorithm requires a partition or a cover of the graph in clusters, this interpolation technique is fused in an innovative way with the GNN clustering methods previously introduced.

\section{Numerical Experiments}

In this section, we provide numerical experiments for signal interpolation using GNNs for community detection. The main MATLAB code for the GBF-PUM method is implemented and available in \cite{cavoretto_2}. After selecting nine unsupervised GNN methods, we compare their performance on two test datasets corresponding to the urban networks of Bologna and Paris \cite{boeing}.

To construct the characteristic matrix for the test graphs, we start from the matrix $\textbf{X} \in \mathbb{R}^{n \times 2}$ containing the two-dimensional node coordinates and multiply it by a random matrix $\textbf{W} \in \mathbb{R}^{2 \times k}$ whose entries are sampled from a normal distribution. We then define 
\begin{equation*}\overline{\textbf{X}} = \cos(\textbf{X}\textbf{W}) \in \mathbb{R}^{n \times k}.
\end{equation*}

If the clustering method requires the number $J$ of communities to be specified in advance, as is the case for all methods except DGCluster, we select the value of $J$ that yields the highest modularity. 
For interpolation, we chose a variational spline as the GBF, uniquely identified from its Fourier transform
\begin{equation*} \label{spline} \hat{\textbf{f}}_{\epsilon \textbf{I}_n + \textbf{L}^{-s}} = \left(\frac{1}{(\epsilon + \lambda_1)^s}, \dots, \frac{1}{(\epsilon + \lambda_n)^s}\right).
\end{equation*}
The chosen parameters for each graph are as follows:
\begin{itemize}
\item For Bologna, $N = 500$, $k = 2000$, $R = 8$, $\epsilon = 0.001$, $s = 2$, $\textbf{x} = \sum_{i = 1}^{10} \textbf{u}_i$;
\item For Paris, $N = 800$, $k = 2000$, $R = 10$, $\epsilon = 0.001$, $s = 2$, $\textbf{x} = \sum_{i = 1}^{10} \textbf{u}_i$.
\end{itemize}

Table \ref{table:res} presents the results for GBF-PUM interpolation on the two test datasets with different GNN models, evaluating the interpolation $\textbf{x}_*$ of signal $\textbf{x}$ using the Relative Maximum Absolute Error (RMAE) and the Relative Root Mean Squared Error (RRMSE), defined as
\begin{equation*}\mbox{RMAE} = \frac{||\textbf{x} - \textbf{x}_*||_{\infty}} {||\textbf{x}||_{\infty}}, \hspace{0.5cm} \mbox{RRMSE} = \frac{||\textbf{x} - \textbf{x}_*||_{2}} {\sqrt{n} \cdot ||\textbf{x}||_{2}}. \end{equation*}

\renewcommand{\arraystretch}{1.5}
\begin{table}[h!]
\centering
\scalebox{0.6}{
\centering
\begin{tabular}
{ |p{1.8cm}|| p{1.5cm} | p{1.5cm}| p{2cm}| p{2cm}| }
 \hline
 Models         & $J$   & MOD     & RMAE   & RRMSE \\
 \hline
 AGC            & 57 & 0.9007 & 4.675e-01     & 4.582e-03\\
 \hline
 AMIL           & 42 & 0.9164 & 1.705e-01     & 1.067e-03\\
 \hline
 AOCD          & 12 & 0.7844 & 1.198e-01     & 6.238e-04\\
 \hline
 CDBNE         & 18 & 0.6115 & 1.263e-01     & 5.943e-04\\
 \hline
 DGCluster     & 7 & 0.8668 & 2.428e-01     & 6.961e-04\\
 \hline
 MAGAE        & 33 & 0.8988 & 1.575e-01     & 1.053e-03\\
 \hline
 MAVGAE      & 42 & 0.8927 & 1.776e-01     & 1.282e-03\\
 \hline
 MGAE          & 36 & 0.9060 & 2.448e-01     & 1.801e-03\\ 
 \hline
 NOCD          & 16 & 0.7542 & 2.662e-01     & 9.109e-04\\
 \hline
\end{tabular}}
\hspace{0.2cm}
\scalebox{0.6}{\centering \begin{tabular}{ |p{1.8cm}|| p{1.5cm} | p{1.5cm}| p{2cm}| p{2cm}| }
 \hline
 Models         & $J$   & MOD     & RMAE    & RRMSE\\
 \hline
 AGC             & 40 & 0.9127 & 5.247e-01     & 9.158e-04\\
 \hline
 AMIL            & 37 & 0.9185 & 7.998e-02 & 1.380e-04\\
 \hline
 AOCD          & 12 & 0.7579 & 5.764e-02 & 1.361e-04\\
 \hline
 CDBNE         & 16 & 0.5923 & 7.560e-02 & 1.012e-04\\
 \hline
 DGCluster     & 9 &  0.8673 & 1.243e-01     & 9.998e-05\\
 \hline
 MAGAE         & 38 & 0.9001 & 6.865e-02 & 1.310e-04\\
 \hline
 MAVGAE       & 40 & 0.9020& 8.619e-02  & 1.201e-04\\
 \hline
 MGAE          & 37 & 0.9086 & 5.425e-02 & 1.334e-04 \\
 \hline
 NOCD          & 17 & 0.7253 & 3.842e-01      & 2.831e-04\\
 \hline
\end{tabular}}
\vspace{0.2cm}

 \caption{\small GBF-PUM interpolation results with communities generated via different GNN methods on the Bologna (left) and Paris (right) graphs.}
\label{table:res}
\end{table}

The results show a good approximation of the test signal for most of the methods, with errors similar or lower than ones obtained using the method presented in \cite{cavoretto_2}. On the vast majority of results, especially for Paris, we can see a large difference in terms of order of magnitudes between the two errors. This means that the difference between the original and the reconstructed signal is low for a vast majority of the nodes, while the error is mostly focused on a few nodes. The optimal number of communities varies substantially across methods, but it does not appear to have a significant impact on the interpolation accuracy.

\section{Conclusions and Future Works}

In this paper, we have investigated the use of GNNs for community detection in graph signal interpolation. The considered GNN-based methods are able to extract meaningful structural information from graphs and generate effective community partitions. 
By incorporating them into a GBF-PUM interpolation framework, where they define the local subdomains used to construct the global approximation, we show with numerical results that this combination can produce accurate signal reconstructions. 

Several directions remain open for future research. First, the need to manually select the number of communities can be computationally demanding, especially when several values of $J$ must be tested to maximize modularity. A natural extension would be to further investigate methods such as DGCluster, which automatically determine the number of clusters. Second, the use of dense characteristic matrices may lead to high computational costs for large graphs. Therefore, another promising direction is the construction of sparse characteristic matrices that preserve the quality of the interpolation while improving the efficiency of the overall procedure.


\section{Acknowledgments}
This work has been supported by the INdAM Research group GNCS, the GFI 2025 Project, and the 2024 Project \lq\lq Numerical Analysis and Modelling\rq\rq\ funded by the Department of Mathematics \lq\lq Giuseppe Peano\rq\rq\ of the University of Torino. This research has been accomplished within the RITA \lq\lq Research ITalian network on Approximation\rq\rq, the UMI Group TAA \lq\lq Approximation Theory and Applications\rq\rq, and the SIMAI Activity Group ANA$\&$A \lq\lq Numerical and Analytical Approximation of Data and Functions with Applications\rq\rq.


%
%
%
%

\end{document}